# On the asymptotics of a cotangent sum related to the Estermann zeta function

George Fikioris

School of Electrical and Computer Engineering

National Technical University of Athens

GR 157-73 Zografou, Athens, Greece

email: gfiki@ece.ntua.gr

*Abstract*

The sum $c_0(1/k) = -\sum_{m=1}^{k-1}(m/k)\cot(m\pi/k)$ is related to the Estermann zeta function. A recent paper computes the first two terms of the large-$k$ asymptotic expansion of $c_0(1/k)$. Using the Poisson summation formula for finite sums, we find three additional terms.

## 1. Introduction and Main Result

For $k = 1, 2, \ldots$, let

$$c_0\left(\frac{1}{k}\right) = -\sum_{m=1}^{k-1} \frac{m}{k} \cot\frac{m\pi}{k} = \sum_{m=1}^{k-1} g_0\left(\frac{m}{k}\right), \qquad (1)$$

where

$$g_0(x) = -x\cot(\pi x). \qquad (2)$$

Our notation $c_0(1/k)$ is the same as that of [1]. As explained in [1], $c_0(1/k)$ is connected to the computation of the zeros of the Estermann zeta function, a problem that is, in turn, related to the Riemann hypothesis. Here, our main result is that $c_0(1/k) - k\ln k/\pi$ can be expanded into an asymptotic power series, with

$$c_0\left(\frac{1}{k}\right) = \frac{1}{\pi}k\ln k - \frac{k}{\pi}\left[\ln(2\pi) - \gamma\right] + \frac{1}{\pi} + \frac{\pi}{36k} - \frac{\pi^3}{5400k^3} + O\left(\frac{1}{k^5}\right), \quad k \to \infty, \quad (3)$$

where $\gamma$ is the Euler-Mascheroni constant. Eqn. (3)—which gives a five-term asymptotic approximation to $c_0(1/k)$— improves upon Theorem 1.3 of [1] (see also [2]), which gives the first two terms. We indicate how to explicitly obtain the next few terms (not given here for brevity) and how to express the remainder of (3) as an integral.

Our derivation is a straightforward application of the Poisson summation formula (PSF) to an auxiliary sum $c(1/k)$. We use the version of the PSF for finite sums; as is well known, this PSF is closely related to the Euler-Maclaurin formula and to the composite trapezoidal rule for quadrature [3]-[6]. The auxiliary sum $c(1/k)$ is associated with a function $g(x)$ that is smooth in $[0,1]$. (Our PSF cannot be applied directly to $c_0(1/k)$ because the associated $g_0(x)$ is infinite when $x = 1$.)

Eqn. (3) gives excellent numerical accuracy, e.g. 8 correct significant digits when $k = 10$; and 3-4 digits even for $k$ as low as 3 (3 is the smallest value of $k$ for which $c_0(1/k)$ is nonzero).

## 2. Derivations and Remarks

For a smooth function $g(x)$ defined in $[0,1]$, the PSF is [3], [4]

$$\sum_{m=1}^{k-1} g\left(\frac{m}{k}\right) = -\frac{1}{2}g(0) - \frac{1}{2}g(1) + k\int_0^1 g(x)dx + 2k\sum_{m=1}^{\infty}\int_0^1 g(x)\cos(2\pi mkx)dx. \quad (4)$$

We apply (4) to

$$g(x) = g_0(x) - \frac{1}{\pi(1-x)} = -x\cot(\pi x) - \frac{1}{\pi(1-x)}, \quad 0 \leq x \leq 1. \quad (5)$$

It is a consequence of (5) that

$$g(0) = -\frac{2}{\pi}, \quad g(1) = -\frac{1}{\pi}, \quad \int_0^1 g(x)dx = -\frac{\ln(2\pi)}{\pi}.$$

Accordingly, (4) gives

$$c\left(\frac{1}{k}\right) = \frac{3}{2\pi} - \frac{k\ln(2\pi)}{\pi} + 2k\sum_{m=1}^{\infty}\int_{0}^{1}g(x)\cos(2\pi mkx)dx, \qquad (6)$$

where

$$c\left(\frac{1}{k}\right) = \sum_{m=1}^{k-1}g\left(\frac{m}{k}\right) = c_0\left(\frac{1}{k}\right) - \frac{k}{\pi}\sum_{m=1}^{k-1}\frac{1}{k-m} = c_0\left(\frac{1}{k}\right) - \frac{k}{\pi}\left[\psi(k)+\gamma\right]. \qquad (7)$$

In (7), the last equality followed from [5]

$$\psi(n+1) = \sum_{m=1}^{n}\frac{1}{m} - \gamma, \quad n = 0,1,2,\ldots.$$

Eqns. (6) and (7) give

$$c_0\left(\frac{1}{k}\right) = \frac{k}{\pi}\left[\psi(k)+\gamma-\ln(2\pi)\right] + \frac{3}{2\pi} + 2k\sum_{m=1}^{\infty}\int_{0}^{1}g(x)\cos(2\pi mkx)dx. \qquad (8)$$

Successive integrations by parts transform the integral in (8) as follows,

$$\int_{0}^{1}g(x)\cos(2\pi mkx)dx = \frac{g'(1)-g'(0)}{(2\pi mk)^2} - \frac{g'''(1)-g'''(0)}{(2\pi mk)^4} - \frac{1}{(2\pi mk)^5}\int_{0}^{1}g^{(5)}(x)\sin(2\pi mkx)dx. \qquad (9)$$

Eqn. (5) implies

$$g'(0) = -\frac{1}{\pi}, \quad g'(1) = \frac{\pi}{3}, \quad g'''(0) = -\frac{6}{\pi}, \quad g'''(1) = \frac{2\pi^3}{15}.$$

Since we also have

$$\sum_{m=1}^{\infty}\frac{1}{m^2} = \frac{\pi^2}{6}, \quad \sum_{m=1}^{\infty}\frac{1}{m^4} = \frac{\pi^4}{90},$$

(8) and (9) yield

$$c_0\left(\frac{1}{k}\right) = \frac{k}{\pi}\left[\psi(k)+\gamma-\ln(2\pi)\right] + \frac{3}{2\pi} + \frac{\pi^2+3}{36\pi k} - \frac{\pi^4+45}{5400\pi k^3} + r(k), \qquad (10)$$

in which the remainder $r(k)$ is given (exactly) by

$$r(k) = \frac{-1}{16\pi^5 k^4} \sum_{m=1}^{\infty} \frac{1}{m^5} \int_0^1 g^{(5)}(x) \sin(2\pi mkx) dx . \qquad (11)$$

One more integration by parts shows that $r(k) = O(1/k^5)$. The desired expansion (3) now follows from (10) and [5]

$$\psi(k) = \ln k - \frac{1}{2k} - \frac{1}{12k^2} + \frac{1}{120k^4} + O\left(\frac{1}{k^6}\right), \quad k \to \infty . \qquad (12)$$

If in (11) we interchange the order of summation and integration and explicitly evaluate the resulting sum, we can obtain the remainder in integral form.

Since $g(x)$ is infinitely differentiable in $[0,1]$, the process of integration by parts can be continued indefinitely. Thus $c_0(1/k) - k\ln k/\pi$ (as well as $c_0(1/k) - k\psi(k)/\pi$) possesses an asymptotic power series. In fact, the process can explicitly give the next few terms of this series, but it soon becomes cumbersome.